\documentclass[12pt]{article}
\usepackage{amssymb,amsmath,amsfonts,amsthm}

    \newtheorem{theorem}{Theorem}[section]
    \newtheorem{lemma}[theorem]{Lemma}

\begin{document}
\title{Perturbation from symmetry and multiplicity of solutions for elliptic problems with
subcritical exponential growth in $\mathbb{R} ^2$}

\author{Cristina Tarsi \thanks{ e-mail: tarsi@mat.unimi.it. \ The
author is member of the research group G.N.A.M.P.A of the Italian
Istituto Nazionale di Alta Matematica (INdAM).} ,\\
Dipartimento di Matematica, Universit\`{a} degli Studi, \\
I-20133 Milano, Italy}
\date{}

\maketitle
 We consider the
following boundary value problem
\[
\left\{
\begin{array}{ll}
 -\Delta u= g(x,u) + f(x,u) %
& \;\;\;\;\;x\in \Omega \\
& \\
u=0 & \;\;\;\;\; x\in \partial \Omega
\end{array}
\right.
\]
where $g(x,-\xi )=-g(x,\xi )$ and $g$ has subcritical exponential
growth in $\mathbb{R} ^2$. Using the method developed by Bolle, we
prove that this problem has infinitely many solutions under
suitable conditions on the growth of $g(u)$ and $f(u)$.

 \medskip

{\bf AMS subject classification:} Primary 35J60, secondary  58E05.

\section{Introduction}

In the last few years, many authors have widely investigated
existence and multiplicity of solutions for semilinear elliptic
problems with Dirichlet boundary conditions by using variational
methods  and topological arguments (see \cite{St1} and references
therein). In particular, the following model
\begin{equation}
\label{model} \left\{
\begin{array}{ll}
 -\Delta u= |u|^{p-2} u + f(x) %
& \;\;\;\;\;x\in \Omega \\
& \\
u=0 & \;\;\;\;\; x\in \partial \Omega
\end{array}
\right.
\end{equation}
 has been extensively studied, where $\Omega$ is an open bounded domain
 of $\mathbb{R} ^N$, $N\geq3$, $f \in L^2 (\Omega )$ and $2<p<2N/N-2$. If $f\equiv 0$,
equation (\ref{model}) possesses natural symmetry, which
guarantees the existence of an unbounded sequence of critical
values for the symmetric functional associated to the problem. On
the contrary, if $f\neq 0$ the problem loses its $\mathbb{Z} _2$
symmetry and a natural question is whether the infinite number of
solutions is preserved or not under perturbation of the odd
equation; a partial answer was independently obtained by
Rabinowitz \cite{Ra}, Bahri-Berestycki \cite{BB} and Struwe
\cite{St2}, who showed in important works that the multiplicity
structure can be maintained also in the perturbed case,
restricting the growth range of the nonlinearity with suitable
bounds depending on $N$. The main idea is to think of the
non-symmetric functional $I$ under study as a perturbation of its
symmetric part $I_0$ and then to estimate how the growth of rate
of the critical levels of $I_0$
is affected by perturbation from symmetry $I-I_0$.\\
More recently, a new type of perturbation from symmetry has been
considered, resulting from second order systems with
non-homogeneous boundary conditions: if $f=0$ in (\ref{model}) but
$u|_{\partial \Omega } =u_0 \neq 0$, the symmetry is again broken
and the perturbation - due to the non-homogeneous boundary
condition - is of higher order. The standard perturbative method
can be applied but yields the result for even smaller range of $p$
values. It was to deal with this type of perturbation that Bolle
\cite{Bo} developed his new approach: this new method deals with
$I$ as the end-point of a continuous path of functionals
$I_{\theta }$, $\theta \in [0,1]$, which starts at the symmetric
functional $I_0$. Bolle's abstract theorem states, roughly
speaking, that the preservation of the min-max critical levels
along the path of functionals $I_{\theta }$ depends only on the
velocities of deformation $\frac{\partial }{\partial \theta
}I_{\theta } (u)$ at the critical points $u$ of $I_{\theta }$.
This fact often allows to obtain better estimates at such points
since they obey certain conservation laws, being solutions of the
corresponding Euler-Lagrange equations. Bolle, Ghoussoub and
Tehrani  \cite{BGT} tested this approach on several other
problems, including the non-homogeneous problem

\begin{eqnarray*}
 \left\{
\begin{array}{ll}
 -\Delta u= |u|^{p-2} u %
& \;\;\;\;\;x\in \Omega \\
& \\
u=u_0 & \;\;\;\;\; x\in \partial \Omega ,
\end{array}
\right.
\end{eqnarray*}

proving the existence of infinitely many solutions for a larger
range, namely for $1<p<(N+1)/(N-1)$. Later, Chambers and Ghoussoub
\cite{CG} have applied Bolle's approach  to establish  a general
multiplicity result for problems with broken symmetry, where the
forcing term $f$ depends also on $u$; they have been able to prove
that the infinite sequence of critical values is preserved if $p$
belongs to a range of values depending also on the growth of
$f(u)$: roughly speaking, the lower is the growth of the
perturbation $f(u)$, the better is the result obtained. \\
In this paper we deal with an analogous of problem (\ref{model})
in dimension $N=2$. Let $\Omega$ be an open bounded subset of
$\mathbb{R}^2$ with smooth boundary $\partial \Omega$; we are
concerned  with existence and multiplicity results for nonlinear
elliptic equations of the type
\begin{equation}
\label{model2} \left\{
\begin{array}{ll}
 -\Delta u= g(x,u) + f(x,u) %
& \;\;\;\;\;x\in \Omega \\
& \\
u=0 & \;\;\;\;\; x\in \partial \Omega
\end{array}
\right.
\end{equation}
where $g(x,-\xi )=-g(x,\xi )$ and $g$ has subcritical growth in
$\mathbb{R} ^2$. When $N=2$, the notion of criticality, that is,
the maximal growth on $u$ which allows to treat problem
(\ref{model2}) variationally, is motivated by the so called
Trudinger-Moser inequality \cite{Tr}, \cite{M}, which says that
for $\alpha \leq 4\pi $
\[
\sup_{\| u\| _{H_0^1} \leq 1} \int_B e^{\alpha u^2}\leq c ( \alpha
)  |B| \leq c( 4\pi ) |B|=C_{TM}|B|
\]
where $\left| B\right| $ denotes the Lebesgue measure of $B$ and
$C_{TM}$ is a constant which does not depend on $u$; hence, the
maximal growth permitted to study problem (\ref{model2})
variationally is of exponential type. Motivated by the Trudinger-
Moser inequality, we say that $g$ has subcritical growth at
$+\infty$ if for all $\alpha >0$
\[
\lim_{\xi \rightarrow +\infty } \frac{|g(\xi )|}{e^{\alpha t^2 } }
=0,
\]
and $g$ has critical growth at $+\infty$ if there exists $\alpha
_0 >0$ such that
\[
\lim_{\xi \rightarrow +\infty } \frac{|g(\xi )|}{e^{\alpha t^2 } }
=0\;\;\;\forall \alpha >\alpha _0; %
\lim_{\xi \rightarrow +\infty } \frac{|g(\xi )|}{e^{\alpha t^2 } }
=+\infty \;\;\;\forall \alpha <\alpha _0.
\]
We will consider only the subcritical case. To our knowledge, the
problem of perturbation from symmetry for equation with
exponential growth in bounded domain of $\mathbb{R} ^2$ has been
approached only by Sugimura \cite{Sug}, who proved that the
infinite number of solutions is preserved if the nonlinear term
has an exponential growth of the kind $e^{\xi ^q}, 0<q<1/2$, and
the forcing term $f=f(x)$ does not depend on $u$. In this paper we
approach the problem (\ref{model2}) using Bolle's method:
following the idea in [CG], we are able to extend the result of
Sugimura to perturbed problem with forcing term depending also on
$u$. As just remarked, maximal growth allowed depends now on the
growth of $f(u)$: in particular we prove the existence of infinite
solutions for (\ref{model2}) if, roughly speaking, $g(u)\sim
e^{\xi ^q}, 0<q\leq 1$ and $f(u)$ satisfies suitable growth's
conditions. Our result includes the one obtained by Sugimura as a
special case: we remark that in this case we are able to include
also the case $q=1/2$, which was not considered in [Sug].\\
Whether the result would still hold for all $q$ up to the
Trudinger-Moser critical exponent $q=2$ is still open.
\section{Preservation of critical levels under deformation of an even functional}
In this section we recall Bolle's method for dealing with problems
with deformation from symmetry (see e.g. \cite{BGT}, \cite{CG}).
Consider two continuous functions $\rho _1$, $\rho _2:[0,1]\times
\mathbb{R} \rightarrow \mathbb{R}$ which are Lipschitz continuous
relative to the second variable. Assume $\rho _1\leq \rho _2$ and
denote by $\psi _1,\psi _2$ the scalar fields associated to them,
defined on $[0,1]\times \mathbb{R}$ by

\[
\left\{
\begin{array}{l}
  \psi _i(0,s)= s \\
   \\
  \frac \partial {\partial \theta} \psi _i(\theta ,s)=\rho
  _i(\theta ,\psi _i(\theta ,s)).
\end{array}
\right.
\]

Note that $\psi _1$ and $\psi _2$ are continuous and that for all
$\theta \in [0,1]$, $\psi _1(\theta, \cdot)$ and $\psi _2(\theta,
\cdot)$ are non decreasing on $\mathbb{R}$; moreover, since $\rho
_1 \leq \rho _2$ we have $\psi _1\leq \psi _2$.
\newline Let $E$ be a Hilbert space
(with scalar product $\langle ,\rangle$ and associated norm $\|
\cdot \|)$ and consider a $\mathcal{C} ^2$ functional $I_0$ on
$E$; let $I$ be a $\mathcal{C} ^2$ functional $:[0,1]\times E
\rightarrow \mathbb{R}$. For subsets $U\subset V$ of $E$, denote
by
\[
c_U(\theta )= \sup_U I_{\theta}, \;\;\;\; %
c_{V,U}(\theta )= \inf_{g\in S_{V,U}} \sup_{g(U)} I_{\theta}
\]

\noindent where
\[
S_{V,U}=\{ g \in \mathcal{C}(V,E):g(v)=v\;\textrm{for}\;%
v\in U\; \textrm{and}\; g(v)=v \; \textrm{for} \; \| v\|
>R, \textrm{for some} \; R>0 \}.
\]
\noindent We make the following assumptions ( we shall use the
abbreviation $I_{\theta}$ for $I(\theta, \cdot )$):
\begin{itemize}
  \item [(H1)] $I$ satisfies a kind of Palais-Smale condition: for
  every sequence $(\theta _n, v_n)$ in $[0,1]\times E$ such that
  $I(\theta _n, v_n)$ is bounded and $\lim_{n\rightarrow \infty}
  \| I'_{\theta_n}(v_n) \| =0$ there is a subsequence converging
  in $[0,1]\times E$. (The limit $(\theta , v)$ then satisfies
  $I'_{\theta }(v)=0$).
  \item [(H2)] For all $b>0$ there is a constant $C_1 (b)$ such
  that:

  \[
  |I_{\theta }(v)|<b\;\;\; \textrm{implies}
  \;\;\;|\frac{\partial I}{\partial \theta} (\theta, v)| \leq
  C_1 (b)(\| I'_{\theta} (v)\| +1)(\| v\| +1).
  \]

  \item [(H3)] For all critical points $v$ of $I_{\theta}$,

  \[
  \rho_1 (\theta , I_{\theta} (v)) \leq \frac{\partial I}{\partial \theta} %
   (\theta, v) \leq \rho_2 (\theta , I_{\theta} (v))
  \]

  \item [(H4)] There are two closed subsets of $E$, $B$ and $A\subset
  B$ such that
    \begin{itemize}
    \item [(i)] $I_0$ has an upper-bound on $B$ and for some
    $\theta_0 \in [0,1]$
    \[
    \lim_{|v|\rightarrow +\infty,v\in B}\sup_{\theta \in
    [0,\theta_0 ]}I_{\theta}(v)=-\infty.
    \]
    \item [(ii)] $c=c_{B,A}(0)>b=c_A(0)$
    \end{itemize}
    \noindent
    In the sequel we will say that $I_0$ \emph{has a min-max configuration}%
    $(c,b)$ if it satisfies hypothesis (H4).
  \item [(H4')] $I_0$ is even and for any finite dimensional
  subspace $W$ of $E$ and any $\theta$ we have
   \[
   \lim_{\| w\| \rightarrow \infty,w\in W} \;\;\; \sup_{\beta \in [0,\theta ]}
   I(\beta ,w)=-\infty.
   \]

\end{itemize}
Observe that we are assuming implicitly in the above definition
that the starting functional $I_0$ satisfies the Palais-Smale
condition $(H_1)$ and $(H_4)$ for $\theta =0$. Set
\[
\overline{\rho}_1 (\theta ,t)=\sup_{\beta \in [0,\theta ]} \rho_1 %
(\beta ,t),\;\;\;\;\;\;\;\;\overline{\rho}_2 (\theta
,t)=\sup_{\beta \in [0,\theta ]} \rho_2 (\beta ,t).
\]
\noindent The main idea of Bolle's result is the following: If one
assumes a min-max critical level for the initial functional $I_0$,
then the deformation velocities $\rho_1$ and $\rho_2$ will
determine whether this critical level persists along the path. We
are now ready to present a reformulation of Bolle's result due to
Chambers and Ghoussoub (see \cite{CG} for further references and
\cite{Bo} for the original result)
\begin{theorem}
\label{Bolle1} Let $I_0$ be a $\mathcal{C}^2$- functional on $E$
with a min-max configuration $c_{B,A}(0)>c_A(0)$, as defined in
$(H4)$; let $\rho_1 \leq \rho_2$ be two velocity fields and
$\psi_1$, $\psi_2$ be the corresponding scalar flows. If
$\psi_1(\theta_0 ,c_{B,A}(0))>\psi_2(\theta_0 ,c_A(0))$ for some
$\theta_0 \in [0,1]$ then for any path of functionals
$I:[0,1]\times E\rightarrow \mathbb{R}$ satisfying $(H1)$, $(H2)$,
$(H3)$ and $(H4)$ the functional $I_{\theta_0}$ has a critical
point at a level $\overline{c}$ such that:
 \[
 \psi_1 (\theta_0 ,c_{B,A}(0))\leq \overline{c} \leq \psi_2 (\theta_0 ,c_{B,A}(0)).
 \]
\end{theorem}

\noindent Note that if $c=c_{B,A}(0)>b=c_A(0)$, as assumed in
$(ii)$ of $(H4)$, it is standard to show that the functional $I_0$
ha a critical point at level $c$: Boll's theorem assures that this
min-max critical level is preserved along any path of functionals
satisfying the above hypotheses if $\psi_1(\theta_0
,c_{B,A}(0))>\psi_2(\theta_0 ,c_A(0))$. \ Assume now that the
Hilbert space is decomposed as $E=\overline{\bigcup
_{k=0}^{\infty} E_k}$ with $E_0=E_-$ being a finite dimensional
subspace and $(E_k)_{k=1}^{\infty}$ is an increasing sequence of
subspaces of $E$ such that $E_k = E_{k-1} \oplus \mathbb{R} e_k$;
let us set

\[
\mathcal{G}=\{g\in \mathcal{C} (E,E): g\;\; \textrm{is odd and for
a fixed}\;\; R>0 \;\; g(v)=v \;\; \textrm{for} \;\; \| v\| \geq R
\}
\]
\noindent and

\begin{equation}
\label{ck}
 c_k = \inf_{g \in \mathcal{G}} \sup_{v \in g(E_k)} I_0(v).
\end{equation}

\noindent In this framework, the following abstract result can be
proved (see \cite{CG}).
\begin{lemma}
\label{Bolle2}
 Let $\rho_1 \leq \rho_2$ be two velocity fields and let $\psi_1
, \psi_2$ be the corresponding scalar flows. Let $I_0$ be an even
$\mathcal{C}^2$ functional on $E$ and consider the levels $c_k$
associated to $I_0$ defined by (\ref{ck}). Then there is $C>0$,
depending only on $\rho_1, \rho_2$ such that for every $k \in
\mathbb{N}$ and every $\theta \in [0,1]$:
\begin{itemize}
\item [(i)] either $\psi_2 (\theta ,c_k)<\psi_1 (\theta ,c_k)$, or
\item [(ii)] $c_{k+1} -c_k \leq C\theta (\overline{\rho}_1 (\theta ,c_{k+1}) +
\overline{\rho}_2 (\theta ,c_k) +1)$.
\end{itemize}
Moreover, in case $(i)$ there exists a level $\ell_k (\theta )$,
only depending on $I_0$ and $\rho_1, \rho_2$, such that for any
good path of functionals $I$ satisfying hypotheses $(H1)$, $(H2)$,
$(H3)$ and $(H4')$ there exists a critical level $\overline{c}_k
(\theta)$ for $I_{\theta}$ with $\psi_2 (\theta ,c_k)< \psi_1
(\theta ,c_{k+1})\leq \overline{c}_k (\theta) \leq \ell_k
(\theta)$.
\end{lemma}
We are now ready to prove the main result of this paper
\begin{theorem}
\label{mainres} Let $\rho_1 \leq \rho_2$ be two velocity fields
and let $\psi_1 , \psi_2$ be the corresponding scalar flows.
Assume now that the Hilbert space is decomposed as
$E=\overline{\bigcup _{k=0}^{\infty} E_k}$ with $E_0=E_-$ being a
finite dimensional subspace and $(E_k)_{k=1}^{\infty}$ is an
increasing sequence of subspaces of $E$ such that $E_k = E_{k-1}
\oplus \mathbb{R} e_k$. Let $I_0$ be an even $\mathcal{C}^2$
functional on $E$ and consider the levels $c_k$ associated to
$I_0$ defined by (\ref{ck}). We have :
\begin{itemize}
\item [(i)] if $\psi_1 (\theta , c_k)\uparrow +\infty$ as %
$k\rightarrow \infty$, then for every $N\in \mathbb{N}$ there
exists a $\theta_N \in (0,1]$, depending only on $I_0$ and
$\rho_1, \rho_2$, such that for any good path of functionals
$I:[0,1]\times E\rightarrow \mathbb{R}$ satisfying $(H1)$, $(H2)$,
$(H3)$ and $(H4')$ the functional $I_{\theta}$ has at least $N$
distinct critical levels, for any $\theta \in [0, \theta_N]$;
\item [(ii)] if $c_k \geq B_1 + B_2 k(\ln k)^{\overline{\beta}}$
where $\overline{\beta}>0, B_1 \in \mathbb{R}, B_2>0$ and if \\
$\overline{\rho}_i (\theta , t)\leq A_1 + A_2
(\ln{(|t|+1)})^{\overline{\alpha}} (\ln \ln (|t|+1))^{-1}$ where
$\overline{\alpha}\geq 0$ and $A_1 ,A_2 \geq 0$, then $I_1$ has an
unbounded sequence of critical levels provided $\overline{\beta}
\geq \overline{\alpha}$.
\end{itemize}
\end{theorem}

\emph{Proof of Theorem \ref{mainres}}. Theorem \ref{mainres} is a
consequence of Lemma \ref{Bolle2}. \medskip \\
$(i)$ Our aim is to prove that for any $N\in \mathbb{N}$ and for
any good path of functional $I_{\theta}$, there is a $\theta _N
\in (0,1]$ such that the functionals $I_{\theta}$ have $N$
distinct critical levels
$d_1(\theta)<d_2(\theta)<...<d_N(\theta)$, for every $\theta \in
[0,\theta _N]$. We obtain the desired sequence of $N$ critical
levels by induction. Let $C>0$ denote the constant appearing in
Lemma \ref{Bolle2}; define
\[
\eta_k= \inf \{ \theta \in [0,1]: c_{k+1}-c_k \leq C\theta
[\overline{\rho}_1 (\theta , c_{k+1})+\overline{\rho}_2(\theta ,
c_k) +1] \} .
\]
Since $\psi_1(\theta, c_k)\uparrow +\infty$ as $k\rightarrow
+\infty$ by assumption, the sequence $c_k$ is unbounded and for
any $M>0$ there is a $K_M>0$ such that
\[
\psi_1(\theta, c_{k+1})>M\;\;\;\textrm{for all} \;\;\theta \in
[0,1] \;\;\textrm{and}\;\; k>K_M;
\]
let us fix $M_1=1$, $k_1=K_{M_1} + 1$ and $\theta_1 = \eta_{k_1}
/2$. By definition of $\theta_1$, for all $\theta \in
[0,\theta_1]$ the alternative $(ii)$ of Lemma \ref{Bolle2} is not
valid; therefore $(i)$ holds, so that for any $\theta \in
[0,\theta_1]$ and any path of functionals $I:[0,1]\times E
\rightarrow \mathbb{R}$ satisfying $(H1)$, $(H2)$, $(H3)$ and
$(H4')$, the functionals $I_{\theta}$ have critical values
$d_1(\theta)$ with
\[
1=M_1<\psi_2 (\theta ,c_{k_1})< \psi_1 (\theta ,c_{k_1+1})< d_1
\theta) \leq \ell_{k_1} (\theta).
\]
Let now $N\in \mathbb{N}$ and suppose that there is a
$\theta_{N-1} \in (0,1]$ such that for any path of functionals
(satisfying hypotheses of Lemma \ref{Bolle2}) the functionals
$I_{\theta}$ have critical values
$d_1(\theta)<d_2(\theta)<...<d_{N-1}(\theta)\leq
\ell_{k_{N-1}}(\theta)$. Let
\[
M_N>\sup_{\theta \in [0,\theta_{N-1}]} \ell_{k_{N-1}}(\theta)
\]
and let $K_N \in \mathbb{N}$ such that
\[
\psi_1(\theta, c_{k+1}) >M_N \;\;\; \textrm{for all} \;\;\; %
k>K_N,\;\; \theta \in [0,1],
\]
which exists by assumption; define $k_N=K_N +1$ and $\theta_N =
\eta_{k_N} /2$. Again, it is clear that for any $\theta \in
[0,\theta_N ]$ and any good family of functionals $I:[0,1]\times E
\rightarrow \mathbb{R}$, the functionals $I_{\theta}$ have
critical values $d_N(\theta)$ with
\[
M_N<\psi_2 (\theta ,c_{k_N})< \psi_1 (\theta ,c_{k_N+1})< d_N
\theta) \leq \ell_{k_N} (\theta);
\]
but, by definition, $M_N>\sup_{\theta \in [0,\theta_{N-1}]}
\ell_{k_{N-1}}(\theta)$; hence, by hypothesis of induction, we can
conclude that the functionals $I_{\theta}$ have $N$ distinct
critical values satisfying
$d_1(\theta)<d_2(\theta)<...<d_N(\theta)$, that is $(i)$.\medskip \\
$(ii)$ Let us suppose by contradiction that the functional $I_1$
has only finitely many critical levels. Since $\rho_i$ is
Lipschitz continuous in the second variable, then there are
$L_i>0$ such that $|\psi_i -s|\leq \theta
(\overline{\rho}_i(\theta ,s)+L_i)$; hence
\begin{equation}
\label{1} \psi_1(1,s)\geq s-\overline{\rho}_1(1 ,s)-L_1%
\geq s-A_3-A_2 \frac{(\ln(|s|+1))^{\overline{\alpha}}}{\ln \ln
(|s|+1)}.
\end{equation}
Since $c_k$ is unbounded, (\ref{1}) implies that also
$\psi_1(c_{k+1})$ is unbounded; therefore, if we suppose that
$I_1$ has only finitely many solutions, the alternative $(i)$ of
Lemma \ref{Bolle2} can not hold, so that $(ii)$ must be true (with
$\theta =1$). Then we have
\begin{eqnarray*}
c_{k+1} -c_k &\leq& C (\overline{\rho}_1 (1,c_{k+1}) +
\overline{\rho}_2 (1 ,c_k) +1) \\
&\leq& C(A_1+A_2 \frac{(\ln(|c_{k+1}|+1))^{\overline{\alpha}}}{\ln \ln (|c_{k+1}|+1)}+A_1+%
A_2 \frac{(\ln(|c_k|+1))^{\overline{\alpha}}}{\ln \ln
(|c_k|+1)}+1)
\end{eqnarray*}
which implies that there is a $k_1>0$ such that
\begin{equation}
\label{2} c_{k+1}\leq c_k + A_4 \frac{
(\ln{c_{k}})^{\overline{\alpha}}}{\ln \ln c_k}\;\;\; \textrm{for
all}\;\;k>k_1.
\end{equation}
Starting from this relation, one can obtain
\begin{equation}
\label{3} c_k<C_1k \frac{(\ln{k} )^{\overline{\alpha}}}{\ln \ln k}
\;\;\;\; \textrm{for all} \;\;\;k>k_1.
\end{equation}
Estimate (\ref{3}) has been proved by Sugimura \cite{Sug}, hence
we will be brief; let us choose $C_1$ such that (\ref{3}) is
verified for $k=k_1$, and
\[
\frac{A_4}{C_1} \left( 1+\overline{\alpha}  +\ln{C_1}
\right)^{\overline{\alpha}} \leq 1;
\]
assume now that for $k>k_1$ (\ref{3}) is valid. Then, by the
choice of $C_1$
\begin{eqnarray*}
c_{k+1}&\leq& c_k + A_4\frac{(\ln{c_{k}})^{\overline{\alpha}}}{\ln \ln c_k} \\
&\leq& C_1k\frac{(\ln{k})^{\overline{\alpha}}}{\ln \ln k}
+A_4\left[ \ln{C_1} +\ln{k} +\overline{\alpha}\ln{\ln{k}} - \ln
\ln \ln k
\right]^{\overline{\alpha}} \\
&\leq& C_1\frac{(\ln{(k+1)})^{\overline{\alpha}}}{\ln \ln
(k+1)}\left[ k+ \frac{A_4}{C_1} \left( \ln{C_1}+ 1+
\overline{\alpha}  \right)^{\overline{\alpha}} \right] \\
&\leq& C_1 (k+1)\frac{(\ln{(k+1)})^{\overline{\alpha}}}{\ln \ln
(k+1)},
\end{eqnarray*}
that is (\ref{3}) for $k+1$. \\
Now, we recall that, by assumption, $c_k \geq B_1 + B_2 k(\log
k)^{\overline{\beta}}$, which is a contradiction under the further
assumption that $\overline{\beta} \geq \overline{\alpha}$.
\section{Perturbation of a symmetric elliptic problem with
exponential growth} The aim of this section is to prove the
existence of infinitely many solutions for the following perturbed
elliptic problem
\begin{equation} \label{p1}
\left\{
\begin{array}{ll}
 -\Delta u= g\left( x,u\right) + f\left( \theta, x,u \right) %
& \;\;\;\;\;x\in \Omega \\
& \\
u=0 & \;\;\;\;\; x\in \partial \Omega
\end{array}
\right.
\end{equation}
where $g(x, \cdot)$ is odd, has exponential growth (as will be
defined) and $f$ is a perturbative term. Let us define $G(x,\xi
)=\int_0^{\xi} g(x,t)dt$ and $F(x,\xi )=\int_0^{\xi} f(x,t)dx$. We
make the following standard assumptions on the symmetric term $g$
(see also \cite{Sug}) :
\begin{itemize}
\item[(g1)] $g\in \mathcal{C} (\overline{\Omega} \times
\mathbb{R}, \mathbb{R})$;
\item [(g2)] there is a constant $A_0>0$ such that
\[
|g(x,\xi )|\leq A_0 e^{\phi (x)}\;\;\;\; \textrm{for}\;\;\; (x,\xi
)\in \overline{\Omega} \times \mathbb{R},
\]
where $\phi:\mathbb{R} \rightarrow \mathbb{R}$ is a function
satisfying $\phi (\xi )\xi ^{-2}\rightarrow 0$ as $|\xi |
\rightarrow +\infty$;
\item [(g3)] there are constants $\mu>0$ and $r_0 \geq 0$ such
that $0<G(x,\xi ) \ln{G(x, \xi )} \leq \mu \xi g(x,\xi )$ for
$x\in \overline{\Omega}$ and $|\xi |\geq r_0$;
\item [(g4)] $g(x,-\xi )=-g(x,\xi )$ for $(x,\xi
)\in \overline{\Omega} \times \mathbb{R}$
\item [(g5)] there exist $0<\alpha_1 \leq \alpha_2<1,A_1,A_2>0$,
and $B_1,B_2>0$ such that
\[
A_1e^{|\xi |^{\alpha_1}}-B_1\leq G(x,\xi )\leq A_2e^{|\xi
|^{\alpha_2}}+B_2 \;\;\;\; \textrm{for} \;\;\;(x,\xi )\in
\overline{\Omega} \times \mathbb{R},
\]
\end{itemize}
and we make the following assumptions on the perturbative term $f$
\begin{itemize}
\item [(f1)] $f\in \mathcal{C} ([0,1]\times \overline{\Omega} \times
\mathbb{R}, \mathbb{R})$ and $f(0,\cdot , \cdot )=0$;
\item [(f2)] there are $0<\beta_1 <\alpha_1$ and $c_1,c_2>0$ such
that
\[
\left| \frac{\partial}{\partial \theta} F(\theta, x, \xi )\right|
\leq c_1 e^{|\xi|^{\beta_1}} +c_2
\]
\item [(f3)] there are $b>0$ and $\varphi (x) \in L^s (\Omega)$ for some $s\geq 1$ such that
\[
f(\theta, x, \xi )\leq b e^{|\xi |^{\beta_1}}+\varphi (x);
\]
\item [(f4)] there are constants $p>0, r>0$ and $d_1,d_2>0$ such
that
\[
\displaystyle \left| \int_{\Omega} \frac{\partial}{\partial
\theta} F(\theta, x,\xi )dx \right| \leq d_1 \| u\|_p^r + d_2
\]
whenever $u$ is a solution of (\ref{p1}).
\end{itemize}

Using the above notation, we have
\begin{theorem}
\label{perturbprobl} Suppose that $g$ satisfies (g1)-(g5), and
suppose that the perturbative term $f$ satisfies (f1)-(f4). Then
we have:
\begin{itemize}
\item [(i)] for every $N\in \mathbb{N}$ there
exists a $\theta_N \in (0,1]$ problem (\ref{p1}) has at least $N$
distinct solutions;
\item [(ii)] if in addition $2/\alpha_2 -2 \geq
 r/\alpha_1$, then problem (\ref{p1}) has an infinite number
of solutions for all $\theta \in [0,1]$.
\end{itemize}
\end{theorem}
The above theorem is applicable to, e.g., the following problem:
\begin{equation}\label{sug}
\left\{
\begin{array}{ll}
  -\Delta u = ue^{|u|^q}+\theta f(x)u^{r-1}& \;\;\;\;\;\textrm{on}\;\; \Omega \\
  & \\
  \;\;\;\;\;\;u=0 & \;\;\;\;\;\textrm{on}\;\; \partial \Omega
\end{array}
\right.
\end{equation}
where
\[
0<q<1,\;\;\;\;\;\; 0\leq r \leq 2-2q
\]
and $f\in \mathcal{L} ^s (\overline{\Omega} )$ for some $s> \frac
r {r-1}$. Theorem $3.1$ includes as a special case the result
obtained by Sugimura, which is now extended also to the exponent $q = 1/2$.  \\
\\
 Let $E=H_0 ^1 (\Omega)$ be the completion of
$\mathcal{C} ^\infty _0 (\Omega)$ with respect to the norm
\[
\| u\| = (\int_{\Omega} |\nabla u|^2 dx)^{1/2},
\]
and define the path of functionals $I:[0,1]\times E \rightarrow
\mathbb{R}$ by
\[
I_{\theta} (u) = \int _{\Omega} (\frac{1}{2} |\nabla u|^2
-G(x,u)-F(\theta, x,u))dx;
\]
then $I_0$ is an even functional and the critical points of
$I_{\theta}$ are solutions of problem (\ref{p1}). It is standard
to verify that $I$ satisfies the Palais Smale condition $(H1)$
(see for example...). Lemma $3.2$ assures that condition $(H2)$ is
satisfied. It is also easy to show that condition $ (H4' )$ is
satisfied in any finite dimensional subspace of $E$, thanks to the
exponential growth of the nonlinear term $G(x,u)$. \newline For
each $k$, denote by $E_k$ the subspace of $E$ spanned by the first
$k$ eigenfunctions of $\triangle$; then, let
\[
\mathcal{G} = \{g \in \mathcal{C} (E,E): g \textrm{ is odd and
}g(x) = x \textrm{ for large } \|x \| \},
\]
and set  $c_k = \inf_{g \in \mathcal{G} } \sup_{g(E_k)} I_0 $. As
before, $c_k$ are critical levels of the even functional $I_0$ .
Lemma $3.3$ below shows that $(H3)$ holds with
\begin{eqnarray*}
\rho _1 (\theta ,s) &=& -C \frac{(\ln (|s|+1))^{\frac{r}{\alpha
_1} }}{\ln \ln (|s|+1)},
\\ \rho _2 (\theta ,s) &=& C \frac{(\ln (|s|+1))^{\frac{r}{\alpha _1}
}}{\ln \ln (|s|+1) };
\end{eqnarray*}
on the other hand, it is shown by Sugimura \cite{Sug} that there
are positive constants $B_1, B_2 $ such that $c_k \geq B_1 k (\ln
k)^{\frac{2}{\alpha _2} -2} -B_2$. Therefore we can apply Theorem
$2.3$ $(ii)$ with $\overline{\alpha } = \frac{r}{\alpha _1}$ and
$\overline{\beta } = \frac{2}{\alpha _2} -2$ to obtain that $I_1$
has an infinite number of solutions when $\frac{2}{\alpha _2} -2
\geq \frac{r}{\alpha _1}$, which is the claim in Theorem
$3.1.(ii)$. Theorem $3.1.(i)$ also follows directly from Theorem
$2.3.(i)$ , since $\psi _1 (\theta ,c_k) \uparrow +\infty$ as
$k\rightarrow \infty $ for the $\rho _i$ defined above.

\begin{lemma}
For all $b>0$ there is a constant $C_1 (b)$ such that \newline
$|I_{\theta }(v)|<b $ implies
\[
|\frac{\partial I}{\partial \theta} (\theta, v)| \leq
  C_1 (b)(\| I'_{\theta} (v)\| +1)(\| v\| +1).
  \]
\end{lemma}
\emph{Proof} Let $b>0$ and suppose that $|I_{\theta }(v)|<b $;
then
\begin{equation}
\label{stimab} \int _{\Omega} e^{|u| ^{{\alpha}_1}} dx \leq \int
_{\Omega} e^{|u| ^{{\alpha}_1}} dx + |\int _{\Omega} F(\theta
,x,u)dx| \leq b + \frac{1}{2} \|u \| ^2
\end{equation}.

Combining $(g3)$ with (\ref{stimab}) yields
\begin{eqnarray*}
-\langle I'_{\theta } (u),u\rangle &=& \int _{\Omega}(- |\nabla
u|^2 +g(x,u)u+f(\theta ,x,u)u)dx \\
&\geq&  \int _{\Omega} |\nabla u|^2 dx + \int
_{\Omega} (g(x,u)u - 4 G(x,u))dx \\
&& + \int _{\Omega} (f(\theta ,x,u)u - 4 F(\theta ,x,u))dx -
4I_{\theta} (u) \\
&\geq& A_4 \|u \| ^2 - A_5 - A_6 \int _{\Omega} e^{|u|
^{{\beta}_1}} dx -4b \\
&\geq& A_7 \|u \| ^2 - A_6 \int _{\Omega} e^{|u| ^{{\beta}_1}} dx
-A_8 \\
&\geq& A_7 \|u \| ^2 - \varepsilon \int _{\Omega} e^{|u|
^{{\alpha}_1}} dx -C_{\varepsilon }
 \\
&\geq& A_9 \|u \| ^2 -A_{10}.
\end{eqnarray*}
Therefore,
\begin{equation}
\label{stima2} \| I'_{\theta} (u) \| \|u \| \geq -\langle
I'_{\theta } (u),u\rangle \geq A_9 \|u \| ^2 -A_{10};
\end{equation}
finally, combining $(f2)$ with (\ref{stima2} ) and (\ref{stimab})
\begin{eqnarray*}
|\frac{\partial}{\partial \theta} I(\theta ,u)| &=& |\int
_{\Omega} \frac{\partial}{\partial \theta} F(\theta ,x,u)dx | \\
&\leq& c_1 \int _{\Omega} e^{|u| ^{{\beta}_1}} dx +c_2 |\Omega |
\\
&\leq& \varepsilon c_1 \int _{\Omega} e^{|u| ^{{\alpha}_1}} dx
+c_{\varepsilon} \\
&\leq& A_{11} \| u \| ^2 +A_{12} \\
&\leq& A_{13} \| I'_{\theta} (u) \| \|u \| +A_{14} \\
&\leq& C (\| I'_{\theta} (u) \| +1) ( \|u \| +1)
\end{eqnarray*}
\\
\begin{lemma}
There exists a constant $C>0$ such that if $u \in H_0^1 (\Omega )$
is a critical point of $I_{\theta}$, then
\[
|\frac{\partial I}{\partial \theta} (\theta, u)| \leq
  C \frac{(\ln (|I_{\theta }(u)| +1))^{\frac{r}{\alpha _1}}}{\ln \ln (|I_{\theta }(u)| +1)}.
  \]
\end{lemma}
\emph{Proof}. Let us suppose that $I'_{\theta }(u) =0$; then,
combining $(g3)$, $(g5)$, $(f3)$ and Young's inequality yields
\begin{eqnarray*}
I_{\theta }(u) &=& I_{\theta }(u) -\frac{1}{2}\langle I'_{\theta
}(u),u\rangle \\
&=& \int _{\Omega } (\frac{1}{2} g(x,u)u -G(x,u))dx + \int
_{\Omega } (\frac{1}{2} f(\theta ,x,u)u -F( \theta ,x,u))dx \\
&\geq& (\frac{1}{2\mu } -\varepsilon ) \int _{\Omega } G(x,u)\ln
G(x,u)dx  -\int _{\Omega }(\frac{c_1}{2}e^{|u|^{\beta _1}} +
|u|e^{|u|^{\beta _1}})dx -C_{\varepsilon } \\
&\geq& C_1 \int _{\Omega }|u|^{\alpha _1}e^{|u|^{\alpha _1}} dx
-C_2
\\
&\geq& C_3 \|u \| _p ^{\alpha _1}e^{\|u \| _p ^{\alpha _1}} -C_4,
\end{eqnarray*}
for any $p\geq 1$, where $C_3 ,C_4$ are positive constants
depending only on $p$ and $\Omega$; therefore,
\begin{equation}
\label{stimacritica} \|u \|_p \leq C_1 (p) \frac{(\ln (|I_{\theta
}(u)| +1)^{\frac{1}{\alpha _1}}}{\ln \ln (|I_{\theta }(u)| +1)} +
C_2 (p).
\end{equation}
Finally, applying $(f4)$ and (\ref{stimacritica}) we obtain
\begin{eqnarray*}
|\frac{\partial }{\partial \theta} I(\theta, u)| &=& |\int
_{\Omega} \frac{\partial }{\partial \theta} F(\theta, u)dx| \\
&\leq& d_1 \|u \|_p^r +d_2 \\
&\leq& C_5 \frac{(\ln (|I_{\theta }(u)| +1)^{\frac{r}{\alpha
_1}}}{\ln \ln (|I_{\theta }(u)| +1)} +C_6,
\end{eqnarray*}
that is our thesis.


\begin{thebibliography}{Pert}
\bibitem[AR]{AR}  A. Ambrosetti and P.H. Rabinowitz, \emph{Dual variational methods in
critical point theory and applications}, J. Funct. Anal. $\mathbf{14%
}$ $(1973),$ $349-381.$

\bibitem[BB]{BB}  A. Bahri and H. Berestycki, \emph{A perturbation method in
critical point theory and applications}, Trans. Amer. Math. Soc. $\mathbf{267%
}$ $(1981),$ $1-32.$

\bibitem[BL]{BL}  A. Bahri and P.-L. Lions, \emph{Morse index of some
min-max critical points,} Comm. Pure Appl. Math. $\mathbf{41}$ $(1988),$ $%
1027-1037.$

\bibitem[Bo]{Bo}  P. Bolle, \emph{On the Bolza problem,}%
J. Diff. Equations $\mathbf{152}$ $(1999),$ $274-288.$

\bibitem[BGT]{BGT}  P. Bolle, N. Ghoussoub and H. Tehrani \emph{The multiplicity %
of solutions in non-homogeneous boundary value problems,} %
Manusc. Math. $\mathbf{101}$ $(2000),$ $325-350.$

\bibitem[CG]{CG}  C. Chambers and N. Ghoussoub, \emph{Deformation
from symmetry and multiplicity of solutions in non-homogeneous problems,}%
Discrete Contin. Dinam. Systems $\mathbf{8}$ $(2002),$ $267-281.$

\bibitem[M]{M}  J. Moser, \emph{A sharp form of an inequality by N. Trudinger%
}, Ind. Univ. Math. J. $\mathbf{20}$ $\left( 1971\right) $,
$1077-1092.$

\bibitem[Ra]{Ra}  P. H. Rabinowitz, \emph{Minimax methods in critical point
theory with applications to differential equations}, Amer. Math.
Soc., Providence, $\left( 1986\right) .$

\bibitem[Ra1]{Ra1}  P. H. Rabinowitz, \emph{Multiple critical points of
perturben symmetric functionals}, Trans. Amer. Math. Soc. $\mathbf{272}$ $%
\left( 1982\right) $, $753-770$.

\bibitem[St1]{St1}  M. Struwe, {\em Variational methods\/}, (Springer-Verlag,
Berlin - Heidelberg - New York, (1990)) .

\bibitem[St2]{St2}  M. Struwe, \emph{Infinitely many critical points for
functionals which are not even and applications to superlinear
boundary value problems}, Manusc. Math. $\mathbf{32}$ $(1980)$,
$335-364.$

\bibitem[Sug]{Sug}  K. Sugimura, \emph{Existence of infinitely many
solutions for a perturbed elliptic equation with exponential
growth}, Nonlin. Anal., Theory, Meth. Appl., $\mathbf{22}$
$(1994)$, $277-293.$

\bibitem[Ta]{Ta}  K. Tanaka,\emph{\ Morse indices at critical points related
to the symmetric mountain pass theorem and applications}, Commun.
partial diff. Equat. $\mathbf{14}$, $\left( 1989\right) $,
$99-128.$

\bibitem[Tr]{Tr}  N. S. Trudinger, \emph{On imbedding into Orlicz spaces and
some applications}, J. Math. Mech. $\mathbf{17}$ $\left( 1967\right) $, $%
473-484.$\
\end{thebibliography}
\end{document}